\documentclass[a4paper,11pt,reqno]{amsart}
\usepackage{graphicx}
\usepackage[latin1]{inputenc}
\usepackage{mathrsfs}
\usepackage{dsfont}
\usepackage{hyperref}
\usepackage{amsmath}
\usepackage{amssymb}
\usepackage{amsthm}
\usepackage{amsfonts}
\usepackage{amstext}
\usepackage{amsopn}
\usepackage{amsxtra}
\usepackage{mathrsfs}
\usepackage{dsfont}
\usepackage{esint}
\usepackage{enumitem}

\newtheorem{thm}{Theorem}
\newtheorem{conjecture}{Conjecture}

\newcommand{\R}{\mathbb{R}}

\newcommand{\C}{\mathbb{C}}

\newcommand{\ii}{\infty}
\newcommand\1{{\ensuremath {\mathds 1} }}
\renewcommand\phi{\varphi}

\renewcommand{\geq}{\geqslant}
\renewcommand{\leq}{\leqslant}

\newcommand{\be}{\begin{equation}}
\newcommand{\ee}{\end{equation}}
\newcommand{\bq}{\begin{equation}}
\newcommand{\eq}{\end{equation}}

\usepackage{color}

\title[An open problem in relativistic Quantum Mechanics]{Which nuclear shape generates the strongest attraction on a relativistic electron?\\ \medskip An open problem in relativistic Quantum Mechanics}

\author[M.J. Esteban]{Maria J. Esteban}
\address{CEREMADE, CNRS, Université Paris-Dauphine, PSL Research University, Place de Lattre de Tassigny, 75016 Paris, France} 
\email{esteban@ceremade.dauphine.fr}

\author[M. Lewin]{Mathieu Lewin}
\address{CEREMADE, CNRS, Université Paris-Dauphine, PSL Research University, Place de Lattre de Tassigny, 75016 Paris, France} 
\email{mathieu.lewin@math.cnrs.fr}

\author[\'E. Séré]{\'Eric Séré}
\address{CEREMADE, Université Paris-Dauphine, PSL Research University, CNRS, Place de Lattre de Tassigny, 75016 Paris, France} 
\email{sere@ceremade.dauphine.fr}

\date{\today}

\DeclareRobustCommand{\SkipTocEntry}[5]{}

\begin{document}

\begin{abstract}
In this article we formulate several conjectures concerning the lowest eigenvalue of a Dirac operator with an external electrostatic potential. The latter describes a relativistic quantum electron moving in the field of some (pointwise or extended) nuclei. The main question we ask is whether the eigenvalue is minimal when the nuclear charge is concentrated at one single point. This well-known property in non-relativistic quantum mechanics has escaped all attempts of proof in the relativistic case. 

\medskip

\end{abstract}

\maketitle


\emph{This article is dedicated to Catriona Byrne on the occasion of her retirement. Her extremely good knowledge of the mathematical community and profession and her kindness made her presence in mathematical events always enjoyable and very useful.}

\section{A conjecture for relativistic electrons}
In this note we describe some conjectures which we recently coined in~\cite{EstLewSer-21a,EstLewSer-21b}, concerning the effect of a nuclear charge on a relativistic electron. We first describe the main conjecture somewhat informally, before we discuss more thoroughly its proper mathematical formulation.
Consider a non-negative finite Borel measure $\mu$ on $\R^3$ and the corresponding linear Schrödinger operator
\begin{equation}
 -\frac{\Delta}{2}-\mu\ast\frac{1}{|x|},
 \label{eq:non-relativistic}
\end{equation}
which describes a non-relativistic electron moving in the Coulomb potential generated by the positive charge distribution $\mu$, in atomic units. The lowest (negative) eigenvalue of this operator is given by the variational principle~\cite{LieLos-01}
\begin{multline}
\lambda_1\left(-\frac{\Delta}{2}-\mu\ast\frac{1}{|x|}\right)\\=\inf_{\substack{\phi\in H^1(\R^3)\\ \int_{\R^3}|\phi|^2=1}}\left\{\frac12\int_{\R^3}|\nabla\phi(x)|^2\,dx-\int_{\R^3}\left(\mu\ast\frac1{|\cdot|}\right)(x)\,|\phi(x)|^2\,dx\right\}. 
\label{eq:Schrodinger_minimum}
\end{multline}
Since this is an infimum over affine functions of $\mu$, we deduce immediately that the eigenvalue is a \emph{concave} function of $\mu$. Therefore, it is minimized, at fixed mass $\mu(\R^3)$, when $\mu$ is proportional to a delta and we have 
\begin{equation}
\lambda_1\left(-\frac{\Delta}{2}-\mu\ast\frac{1}{|x|}\right)\geq \lambda_1\left(-\frac{\Delta}{2}-\frac{\mu(\R^3)}{|x|}\right)=-\frac{\mu(\R^3)^2}{2}
 \label{eq:estim_Schrodinger}
\end{equation}
for every $\mu\geq0$. The interpretation is that the lowest possible electronic energy is reached by taking the most concentrated charge distribution, at fixed total charge $\mu(\R^3)$. In fact, in~\cite{LieSim-78,Lieb-82} it is proved that the eigenvalue decreases when $\mu$ is deformed using an arbitrary contraction, for instance a dilation $\alpha^3\mu(\alpha\cdot)$ with $\alpha\geq1$. This was generalized to molecular systems in~\cite{LieSim-78,HofOst-80,Lieb-82}, where it is proved that the electronic part of the ground state energy decreases when all the distances between the nuclei are decreased.

\medskip

Relativistic effects play an important role in the description of quantum electrons in molecules containing heavy nuclei, even for not so large values of the nuclear charge. A proper description of such systems is based on the Dirac operator~\cite{Thaller,EstLewSer-08}. This is a first-order differential operator which has very different properties compared to its non-relativistic counterpart $-\Delta/2$ in~\eqref{eq:non-relativistic}. For instance the spectrum of the free Dirac operator is not semi-bounded which prevents from giving an unambiguous definition of a  ``ground state'' and turns out to be related to the existence of the positron~\cite{EstLewSer-08}. In addition, because of its scaling properties, the Dirac operator has a critical behavior with respect to the Coulomb potential $1/|x|$ which gives a bound $Z\leq 137$ on the highest possible charge of atoms in the periodic table, for point nuclei.

In atomic units for which $m=c=\hbar=1$, the free Dirac operator $D_0$ can be written as
\begin{equation}
D_0\ = -i\boldsymbol{\alpha}\cdot\boldsymbol{\nabla} + \beta = \ - i\sum^3_{k=1} \alpha_k \partial_{x_k} \ + \ {\bf \beta},
\label{eq:def_Dirac}
\end{equation}
where $\alpha_1$, $\alpha_2$, $\alpha_3$ and $\beta$ are $4\times4$ Hermitian matrices which  satisfy the following anticommutation relations:
\begin{equation*} \label{CAR}
\left\lbrace
\begin{array}{rcl}
 {\alpha}_k
{\alpha}_\ell + {\alpha}_\ell
{\alpha}_k  & = &  2\,\delta_{k\ell}\,\1,\\
 {\alpha}_k {\beta} + {\beta} {\alpha}_k
& = & 0,\\
\beta^2 & = & \1.
\end{array} \right. 
\end{equation*}
The usual representation in $2\times 2$ blocks is given by 
$$ \beta=\left( \begin{matrix} I_2 & 0 \\ 0 & -I_2 \\ \end{matrix} \right),\quad \; \alpha_k=\left( \begin{matrix}
0 &\sigma_k \\ \sigma_k &0 \\ \end{matrix}\right)  \qquad (k=1, 2, 3)\,,
$$
where the Pauli matrices are defined as
\begin{equation}
 \sigma _1=\left( \begin{matrix} 0 & 1
\\ 1 & 0 \\ \end{matrix} \right),\quad  \sigma_2=\left( \begin{matrix} 0 & -i \\
i & 0 \\  \end{matrix}\right),\quad  \sigma_3=\left( 
\begin{matrix} 1 & 0\\  0 &-1\\  \end{matrix}\right) \, .
 \label{eq:Pauli}
\end{equation}
The operator $D_0$ is self-adjoint on the domain $H^1(\R^3,\C^4)$ in the Hilbert space $L^2(\R^3,\C^4)$ and its spectrum is $\sigma(D_0)=(-\ii,-1]\cup[1,\ii)$~\cite{Thaller}.  Moreover, $(D_0)^2= -\Delta +1$.

A relativistic electron in the presence of the nuclear charge $\mu$ is described by the Dirac-Coulomb operator
\begin{equation}
D_0-\mu\ast\frac1{|x|}
\label{eq:def_D0_mu}
\end{equation}
in place of the non-relativistic operator~\eqref{eq:non-relativistic}. In our units $\mu$ represents the nuclear charge multiplied by the fine-structure constant $\alpha\simeq1/137$. We defer the precise definition of the Dirac-Coulomb  operator to the next section. Eigenvalues in the gap $(-1,1)$ physically correspond to stationary states of the relativistic electron. Therefore it seems natural to expect that the lowest eigenvalue in $(-1,1)$ will again be minimized for the Dirac measure $\mu(\R^3)\delta_0$, like in the Schrödinger case~\eqref{eq:estim_Schrodinger}. This is the conjecture which we recently made in~\cite{EstLewSer-21a,EstLewSer-21b}.

\begin{conjecture}[General charges~\cite{EstLewSer-21a,EstLewSer-21b}]\label{conjecture_mu}
For any non-negative Borel measure $\mu$ such that $\mu(\R^3)\leq 1$, the lowest eigenvalue in the gap $(-1,1)$ satisfies
\begin{equation}
\lambda_1\left(D_0-\mu\ast\frac1{|x|}\right)\geq\lambda_1\left(D_0-\frac{\mu(\R^3)}{|x|}\right)=\sqrt{1-\mu(\R^3)^2}.
\label{eq:conjecture_eigenvalue_intro}
\end{equation}
\end{conjecture}

In relativistic quantum chemistry one often relies on extended nuclear charges, hence the interest of looking at any possible $\mu$. If we restrict our attention to pointwise nuclei, then we have $\mu=\sum_m\theta_m\delta_{R_m}$ and the conjecture becomes

\begin{conjecture}[Multi-center potentials~\cite{EstLewSer-21a,EstLewSer-21b}]\label{conjecture_multicenter}
We have 
\begin{equation}
\lambda_1\left(D_0-\sum_{m=1}^M\frac{\theta_m}{|x-R_m|}\right)\geq \lambda_1\left(D_0-\frac{\sum_{m=1}^M\theta_m}{|x|}\right)=\sqrt{1-\left(\sum_{m=1}^M\theta_m\right)^2}
\end{equation}
for all $M\geq2$, all $R_1,...,R_M\in\R^3$ and all $\theta_m\geq0$ so that $\sum_{m=1}^M\theta_m\leq1$,
\end{conjecture}

Since any $\mu$ can be approximated by a combination of Dirac deltas for the narrow topology, Conjecture~\ref{conjecture_multicenter} is equivalent to Conjecture~\ref{conjecture_mu}. Indeed $\lambda_1$ is continuous for this topology~\cite[Lemma 12]{EstLewSer-21b}.

The case $M=2$ was conjectured by Klaus in~\cite[p.~478]{Klaus-80b} and by Briet-Hogreve in~\cite[Sec.~2.4]{BriHog-03}. Numerical simulations from~\cite{ArtSurIndPluSto-10,McConnell-13} seem to confirm the conjecture for $M=2$, even for large values of the nuclear charges. In~\cite{EstLewSer-21a,EstLewSer-21b} and here we make the stronger conjecture that the same holds for any $M$. Note that the numerical simulations seem to indicate that $\lambda_1$ decreases when the Euclidean distance between nuclear charges is decreased, a property proved by Lieb and Simon \cite{LieSim-78,Lieb-82} in the non-relativistic case. This leads to a third conjecture:

\begin{conjecture}[Monotonicity]\label{conjecture_monotonicity_mu}
Let $\mu$ be a non-negative Borel measure such that $\mu(\R^3)\leq 1$ and let $f:\R^3\to\R^3$ be a contraction for the Euclidean norm of $\R^3$. Then, denoting by $f_*\mu$ the pushforward of $\mu$ by $f$, we have 
\begin{equation}
\lambda_1\left(D_0-\mu\ast\frac1{|x|}\right)\geq\lambda_1\left(D_0-(f_*\mu)\ast\frac1{|x|}\right).
\label{eq:conjecture_eigenvalue_contraction}
\end{equation}
\end{conjecture}

Conjecture ~\ref{conjecture_mu} is a special case of Conjecture ~\ref{conjecture_monotonicity_mu}, as can be seen by taking $f=0$. In this note we only discuss Conjecture ~\ref{conjecture_mu}, which is already far from obvious. The main difficulty is that the lowest Dirac eigenvalue in the gap $(-1,1)$ is not given by a minimum like in~\eqref{eq:Schrodinger_minimum}. In fact, as quickly explained below, it is given by a min-max formula~\cite{GriSie-99, DolEstSer-00, SchSolTok-20,EstLewSer-21a}. Unfortunately, it does not seem easy to derive a concavity property of $\lambda_1(D_0-\mu\ast|x|^{-1})$ from this variational characterization, and this prevents us from using the same argument as in the nonrelativistic case. However the min-max formula implies that $\lambda_1(D_0+V)$ is monotone in $V$, so that Conjecture~\ref{conjecture_mu} holds true if one restricts it to \emph{radially symmetric measures}~$\mu$. Indeed, for such measures we have the \emph{pointwise} bound
$$\left(\mu\ast\frac1{|\cdot|}\right)(x)\leq \frac{\mu(\R^3)}{|x|},$$
by Newton's theorem~\cite{LieLos-01} and~\eqref{eq:conjecture_eigenvalue_intro} follows. If one only considers radial contractions $f$, Conjecture~\ref{conjecture_monotonicity_mu} is also true for radially symmetric measures~$\mu$. No other case seems to have been proved in the literature. 

In the next section we discuss the proper definition of the Dirac operator $D_0-\mu\ast|x|^{-1}$ in~\eqref{eq:def_D0_mu} and the exact meaning of the ``lowest eigenvalue in the gap'' $\lambda_1(D_0-\mu\ast|x|^{-1})$ appearing in the conjecture.

\section{Dirac operator with external charges}

\subsection{Self-adjointness}
For Coulomb-like potentials $V$, it is not an easy task to define $D_0+V$ as a self-adjoint operator. The reason is that $1/|x|$ has the same homogeneity as the differential part $\boldsymbol{\alpha}\cdot\boldsymbol{\nabla}$ of the free Dirac operator. In the pure Coulomb case $\mu=\nu\delta_0$, everything is explicit. The operator $D_0-\nu|x|^{-1}$ has a unique self-adjoint realization for $\nu\leq\sqrt{3}/2$ and infinitely many for $\nu>\sqrt{3}/2$. For $\nu\in(\sqrt{3}/2,1]$ one self-adjoint extension is special, with the corresponding eigenfunctions being the least singular at the origin. It is called the ``distinguished" extension. For $\nu>1$ all the self-adjoint realizations look the same, with eigenfunctions having similar oscillations near the origin~\cite{Hogreve-13}. For $\nu\in[0,1]$ it is known that the lowest eigenvalue of the distinguished extension in the gap $(-1,1)$ equals $\sqrt{1-\nu^2}$ and therefore remains positive. The formula for this eigenvalue was already used on the right side of~\eqref{eq:conjecture_eigenvalue_intro}. 

Many works have been devoted to the case of a general Coulomb-type potential $V$ since the 70s~\cite{Schmincke-72b,Wust-73,Wust-75,Wust-77,Nenciu-76,Nenciu-77,KlaWus-78,Klaus-80b,Kato-83}. Various methods were introduced to prove that there also exists a unique ``distinguished'' self-adjoint extension. The results typically cover any potential $V$ satisfying the pointwise inequality 
$$0\geq V(x)\geq -\frac{\nu}{|x|},\qquad \nu\in (0, 1).$$ 
In this case, ``distinguished'' can have several possible meanings, which were all eventually shown to be equivalent. One requirement was that the domain of the operator be a subspace of $H^{1/2}( \R^3, \C^4)$, so that the energy is well defined. Another natural property was that the operator is the norm-resolvent limit of the Dirac operator with a regularized potential. Using a quite different approach Esteban and Loss proved more recently in~\cite{EstLos-07,EstLos-08} that a distinguished self-adjoint extension could also be defined in the critical case $\nu=1$. 

For small values of $\nu$, the domain of self-adjointness is just the Sobolev space $H^1(\R^3,\C^4)$ but for larger values of $\nu$, the domain was not explicit in most of the above-cited works. The recent articles~\cite{EstLewSer-19,SchSolTok-20} contain a more detailed analysis of the domain. 

In~\cite{EstLewSer-21a} all the previous works were generalized to cover the case of potentials $V= -\mu\ast|x|^{-1}$. The existence of a ``distinguished'' extension was shown under the sole assumption that $\mu$ is a non-negative finite measure which has no atom of mass larger than or equal to $1$. This gave a clear definition to the operator $D_0-\mu\ast|x|^{-1}$ in~\eqref{eq:def_D0_mu}, describing one electron in the presence of a nuclear charge $\mu$.

\subsection{Dirac eigenvalues in the gap}
Once the operators have been well defined, the next question is how to find and characterize the stationary states, that is, the eigenvalues in the spectral gap $(-1, 1)$. This has also attracted a lot of attention in spectral theory and mathematical physics in the last two decades~\cite{GriSie-99, DolEstSer-00, DolEstSer-00a, DolEstSer-03, DolEstSer-06, MorMul-15,Muller-16,EstLewSer-19,SchSolTok-20}. We are not going to state the precise result here, but the conclusion is that one can characterize the eigenvalues in the spectral gap using non-standard min-max variational methods. Potentials of the form $V=-\mu\ast|x|^{-1}$ were not covered by most of the existing results but they were handled in~\cite{EstLewSer-21a}, following the method in~\cite{DolEstSer-00,EstLewSer-19,SchSolTok-20}. 

Let us emphasize that there is some difficulty in defining what it means to be the ``lowest eigenvalue in the gap $(-1,1)$'', as in our two Conjectures~\ref{conjecture_mu}--\ref{conjecture_monotonicity_mu}. If we have a well-behaved (e.g.~bounded) negative potential $V$, then the eigenvalues of $D_0+tV$ will be close to $1$ for small $t>0$ and will all decrease when $t$ is increased. The lowest eigenvalue will eventually touch the lower spectrum at $-1$, at a certain finite value of $t$, and dissolve in the continuum. Then the second eigenvalue in the gap becomes the lowest one. We do not wish to look at these pathological discontinuities and want to be sure that the lowest eigenvalue remains so for all $t\leq1$.

In fact, should our Conjectures~\ref{conjecture_mu} and~\ref{conjecture_multicenter} hold true, they would imply that
$$\lambda_1\big(D_0-t\mu\ast|x|^{-1}\big)\geq0,\qquad \forall t\in(0,1).$$
In particular, when we turn on the potential $V=-\mu\ast|x|^{-1}$ by means of the parameter $t$, the lowest eigenvalue will always be non-negative and there will be no spectrum in the lower half of the gap $(-1,0)$. No eigenvalue will dive into the negative continuum, which justifies considering the lowest one. 

Since we do not know how to prove the conjecture, a natural first step was to investigate which measure $\mu$ can have eigenvalues approaching the negative threshold $-1$. In~\cite{EstLewSer-21b}, we defined a critical charge  $\nu_1$ as the largest positive number for which 
$$\lambda_1(D_0-\mu\ast|x|^{-1}) > -1 \;\mbox{for all}\; 0<\mu(\R^3)<\nu_1\,.$$
For measures with $\mu(\R^3)<\nu_1$ there is thus no ambiguity of what it means to be the ``lowest eigenvalue''. Our Conjectures~\ref{conjecture_mu} and~\ref{conjecture_multicenter} contain the statement that $\nu_1=1$. The following was shown in~\cite{EstLewSer-21b}.

\begin{thm}[The critical charge $\nu_1$~\cite{EstLewSer-21b}]
The critical number $\nu_1$ satisfies
\begin{equation}
0.9\simeq\frac2{\pi/2+2/\pi}\leq \nu_1\leq1. 
 \label{eq:Tix}
\end{equation}
It is also the best constant in the Hardy-type inequality
\begin{equation}
\boxed{  \int_{\R^3}\frac{|\sigma\cdot\nabla\phi|^2}{\mu\ast|x|^{-1}}\,dx\geq\frac{\nu_1^2}{\mu(\R^3)^2}\int_{\R^3}\left(\mu\ast\frac1{|x|}\right)|\phi|^2\,dx} 
 \label{eq:Hardy_intro}
\end{equation}
for every $\phi\in C^\ii_c(\R^3,\C^2)$ and every finite non-negative measure $\mu\geq0$, where $\sigma_1,\sigma_2,\sigma_3$ are the $2\times2$ Pauli matrices defined above in~\eqref{eq:Pauli}.  
\end{thm}

The estimate~\eqref{eq:Tix} was proved using an inequality due to Tix~\cite{Tix-98}, whereas the link with the Hardy inequality~\eqref{eq:Hardy_intro} comes from the variational characterization of the first eigenvalue. Such inequalities have played an important role in the study of Dirac operators~\cite{DolEstSer-00,DolEstLosVeg-04,DolEstDuoVeg-07,ArrDuoVeg-13,CasPizVeg-20}.

\section{Two results from~\cite{EstLewSer-21a,EstLewSer-21b}}

In this last section we mention two results from~\cite{EstLewSer-21a,EstLewSer-21b} which are related to our Conjectures~\ref{conjecture_mu}--\ref{conjecture_monotonicity_mu}. 

\subsection{Existence of an optimal measure $\mu$}

Even if we do not know that concentrating all the mass at one point gives the lowest eigenvalue, we could at least prove that there exists an optimizer $\mu$ for a fixed mass $\mu(\R^3)=\nu<\nu_1$ and that it has a very small support. 

\begin{thm}[Existence of an optimal measure~\cite{EstLewSer-21b}]
For any $\nu\in[0,\nu_1)$, there exists a positive Borel measure $\mu_\nu$ with $\mu_\nu(\R^3)=\nu$ so that 
$$\lambda_1\left(D_0-\mu_\nu\ast\frac1{|x|}\right)=\min_{\substack{\mu\ :\\ \mu(\R^3)=\nu}}\lambda_1\big(D_0-\mu\ast|x|^{-1}\big).$$
The support of any such minimiser $\mu_\nu$ is a compact set of zero Lebesgue measure.
\end{thm}

The theorem is proved in~\cite{EstLewSer-21b} by a rather delicate adaptation of techniques from nonlinear analysis to the context of Dirac operators. The first eigenvalue is a highly nonlinear function of the measure $\mu$, even if the operator only depends linearly on $\mu$. The main ``enemy" is the action of the non-compact group of space translations, which is controlled using Lions' concentration-compactness method~\cite{Lions-84,Lions-84b}. The main difficulty was to prove that the problem is locally compact under the assumption that $0\leq\nu<\nu_1$ and this is another reason why the critical mass $\nu_1$ plays a central role. In spirit, the local compactness holds true because the eigenvalue cannot dive into the lower continuous spectrum by definition of $\nu_1$. But the actual proof is rather involved and relies on variational arguments using the min-max characterization of the first eigenvalue. That the support has zero Lebesgue measure was shown in~\cite{EstLewSer-21b} by means of a unique continuation principle for Dirac operators, which extends famous results in the Schrödinger case~\cite{JerKen-85,Stein-85}.

\subsection{The potential energy surface}
In quantum chemistry one is interested in the \emph{potential energy surface} which, by definition, is the graph of the first eigenvalue of the multi-center Dirac-Coulomb operator, seen as a function of the locations of the nuclei, including the nuclear repulsion:
$$(R_1,...,R_M)\mapsto \lambda_1\left(D_0-\sum_{m=1}^M\frac{\theta_m}{|x-R_m|}\right)+\sum_{1\leq m<\ell\leq M}\frac{\theta_m\theta_\ell}{|R_m-R_\ell|}.$$
For the case $M=2$ the properties of the above function were analyzed in~\cite{Klaus-80b,HarKla-83,BriHog-03} in the case of subcritical singularities with charge $\theta_m<1$. In~\cite{EstLewSer-21a} we extended these results to cover the case $M>2$ and also to  include the critical case of nuclear charge equal to $1$. We proved the following

\begin{thm}[The potential energy surface~\cite{EstLewSer-21a}]\label{Thm-continuity}
Let $0< \theta_1,...,\theta_M\leq 1$.

\medskip

\noindent $(i)$ The map $(R_1,...,R_M)\mapsto \lambda_1\big(D_0-\sum_{m=1}^M\theta_m|x-R_m|^{-1}\big)$ is  continuous on the open set 
\begin{multline*}
\Omega=\bigg\{(R_1,...,R_M)\in(\R^3)^M\ :\ R_m\neq R_\ell\ \text{for $m\neq \ell$}\\
\lambda_1\left(D_0-\sum_{m=1}^M\frac{\theta_m}{|x-R_m|}\right)>-1\bigg\}.
\end{multline*}


\noindent $(ii)$ Moreover,
\begin{equation}
\lim_{\min_{k\neq \ell}|R_k-R_\ell|\to\ii}\lambda_1\left(D_0-\sum_{m=1}^M\frac{\theta_m}{|x-R_m|}\right)=\sqrt{1-\max_m\theta_m^2}.
\end{equation} 

\medskip

\noindent $(iii)$ If in addition $\sum_{m=1}^M\theta_m< \nu_1$ then 
\begin{equation}
\lim_{\max_{k\neq \ell}|R_k-R_\ell|\to0}\lambda_1\left(D_0-\sum_{m=1}^M\frac{\theta_m}{|x-R_m|}\right)=\sqrt{1-\left(\sum_{m=1}^M\theta_m\right)^2}.
\label{eq:limit_merged}
\end{equation}
\end{thm}

By $(ii)$ we see that Conjecture~\ref{conjecture_multicenter} is valid when the nuclei are infinitely far apart. On the other hand, $(iii)$ says that the lowest eigenvalue is continuous when all the nuclei are merged to one point. Conjecture~\ref{conjecture_multicenter} says that the limit~\eqref{eq:limit_merged} should be from above and it would be interesting to try to prove the conjecture when the nuclei are very close to each other. The limit~\eqref{eq:limit_merged} was also stated for $M=2$ and $\nu_1=\nu_2<1/2$ in~\cite{BriHog-03} but we could not fill all the details of the argument of the proof.

\bigskip

The properties of Dirac-Coulomb operators are fascinating and much more involved than the non-relativistic Schrödinger case. Many tools (such as min-max methods) have been developed to better deal with Dirac operators. Our Conjectures~\ref{conjecture_mu}, ~\ref{conjecture_multicenter} and ~\ref{conjecture_monotonicity_mu} are strongly supported by numerical results in the physics and chemistry literature, but their proof will probably require introducing new techniques. 

\bigskip

\noindent{\textbf{Acknowledgement.}} This project has received funding from the European Research Council (ERC) under the European Union's Horizon 2020 research and innovation programme (grant agreement MDFT No 725528 of M.L.), and from the Agence Nationale de la Recherche (grant agreement molQED).

\bibliographystyle{siam}
\bibliography{biblio}

\end{document}